# The Hilbert Zonotope and a Polynomial Time Algorithm for Universal Gröbner Bases

Eric Babson     Shmuel Onn *     Rekha Thomas †

**Abstract**

We provide a polynomial time algorithm for computing the *universal Gröbner basis* of any polynomial ideal having a finite set of common zeros in fixed number of variables. One ingredient of our algorithm is an effective construction of the state polyhedron of any member of the Hilbert scheme $\text{Hilb}_n^d$ of $n$-long $d$-variate ideals, enabled by introducing the *Hilbert zonotope* $\mathcal{H}_n^d$ and showing that it simultaneously refines all state polyhedra of ideals on $\text{Hilb}_n^d$.

## 1 Introduction

The *universal Gröbner basis* of an ideal $I$ in the algebra $\mathbb{F}[x] := \mathbb{F}[x_1, \ldots, x_d]$ of $d$-variate polynomials over a field is the minimal set $\mathcal{U}(I)$ which is simultaneously a Gröbner basis for $I$ under every monomial order. A finite universal Gröbner basis always exists and in a sense is the ultimate generating set of $I$ for algorithmic purposes. In particular, for ideals having a finite set of common zeros (variety) over the algebraic closure of $\mathbb{F}$, a universal Gröbner basis reduces the problem of computing the zero set to the problem of finding roots of $d$ univariate polynomials. For instance, consider the system $P = \{x_1^2 - x_2, x_2^2 - 7x_2 + 6x_1, x_1x_2 - 3x_2 + 2x_1\}$ of bivariate polynomials; the universal Gröbner basis of the ideal $I := \text{ideal}(P)$ is

$$\mathcal{U}(I) \;=\; P \cup \{x_1^3 - 3x_1^2 + 2x_1, \; x_2^3 - 5x_2^2 + 4x_2, \; x_1 + \frac{1}{6}x_2^2 - \frac{7}{6}x_2, \; x_2 - x_1^2\}$$

and contains univariate polynomials in each of $x_1$ and $x_2$. Finding the roots of these two polynomials, we conclude that the set of zeros of $P$ satisfies $\text{var}(P) = \text{var}(I) \subseteq \{0,1,2\} \times \{0,1,4\}$. Substituting back to $P$, we find that $\text{var}(P) = \{(0,0), (1,1), (2,4)\}$, thereby solving the system.

The *length* of an ideal $I$ in $\mathbb{F}[x]$ is the $\mathbb{F}$-dimension of the quotient $\mathbb{F}[x]/I$ and is finite if and only if the set of common zeros of $I$ over the algebraic closure of $\mathbb{F}$ is finite. Let $\text{Hilb}_n^d$ be the set of $n$-long ideals in $\mathbb{F}[x] = \mathbb{F}[x_1, \ldots, x_d]$; it can be embedded as an algebraic variety in a higher dimensional space and is referred to as the *Hilbert scheme* of $n$-long $d$-variate ideals. One of the goals of this article is to provide a polynomial time algorithm for computing the universal Gröbner basis of any ideal on the Hilbert scheme (see Section 4 for the complete formulation):

*Research partially supported by a grant from ISF - the Israel Science Foundation.
†Research partially supported by DMS grant 0100141.

**Theorem 4.2** Fix $d$. Then there is a polynomial time algorithm that computes the universal Gröbner basis $\mathcal{U}(I)$ of any ideal $I \in \text{Hilb}_n^d$ using $O(n^{2d+1} (\log n)^{(2d-1)(d-1)})$ arithmetic operations.

The computational complexity is measured in terms of the number of arithmetic operations over the underlying field $\mathbb{F}$. Over the field of rational numbers, the algorithm is (strongly) polynomial time in the Turing computation model, but we do not dwell on the details here.

One ingredient of our algorithm is an effective unified construction, for any ideal $I \in \text{Hilb}_n^d$, of its *state polyhedron* $\mathcal{S}(I)$ whose vertices bijectively index the reduced Gröbner bases of $I$. This is done in Section 2, where we introduce the *basis polytope* $\mathcal{B}(I)$ of any ideal $I \in \text{Hilb}_n^d$ and establish the following description of its state polyhedron (see Section 2 for the complete statement):

**Theorem 2.4** The state polyhedron of $I \in \text{Hilb}_n^d$ is provided by $\mathcal{S}(I) := B(I) + \mathbb{R}_+^d$.

As a corollary, we obtain the following polynomial upper bounds on the number of reduced Gröbner bases and the size of the universal Gröbner basis of any ideal on the Hilbert scheme:

**Corollary 2.5** For every fixed $d$, the following hold for any $n$-long $d$-variate ideal $I \in \text{Hilb}_n^d$:

(1) the number of distinct reduced Gröbner bases of $I$ is $O(n^{2d\frac{d-1}{d+1}})$;

(2) the number of elements in the universal Gröbner basis $\mathcal{U}(I)$ is $O(n^{2d-3+\frac{3d-1}{d(d+1)}})$.

The cardinality of the set defining the basis polytope $\mathcal{B}(I)$ of $I \in \text{Hilb}_n^d$ is typically exponential in $n$ even for fixed $d = 2$ and so Theorem 2.4 does not lead directly to an efficient algorithm for constructing the state polyhedron. We overcome this difficulty by introducing, in Section 3, the *Hilbert zonotope* $\mathcal{H}_n^d$. Proving that $\mathcal{B}(I)$ is a projection of a suitable matroid polytope, we show that $\mathcal{H}_n^d$ is *universal* for the Hilbert Scheme in the following sense:

**Theorem 3.5** The Hilbert zonotope $\mathcal{H}_n^d$ refines the state polyhedron $\mathcal{S}(I)$ of every $I \in \text{Hilb}_n^d$.

Using Theorems 2.4 and 3.5 we are able, in Section 4, to provide the aforementioned polynomial time algorithm for constructing the state polyhedron and universal Gröbner basis $\mathcal{U}(I)$ of any $I \in \text{Hilb}_n^d$. In particular, our results apply for the vanishing ideal of *any* point configuration, extending earlier results of [10] for the generic case, and for lattice ideals studied earlier in [11].

In Section 5 we interpret some of the notions and demonstrate some of the results discussed herein for the special classes of vanishing ideals of point configurations and of lattice ideals, the latter having some consequences for the so-called "group relaxation" of integer programming.

We conclude with a brief discussion, in Section 6, of the embedding of the Hilbert scheme $\text{Hilb}_n^d$ into the Grassmanian of $n$-dimensional subspaces of a vector space of dimension $O(n (\log n)^{d-1})$.

## 2   The basis polytope and the state polyhedron

A *staircase* is a set $\lambda \subseteq \mathbb{N}^d$ of nonnegative integer vectors such that $u \leq v \in \lambda$ (coordinatewise) implies $u \in \lambda$. Let $\binom{\mathbb{N}^d}{n}_{stair}$ denote the finite set of $n$-element staircases in $\mathbb{N}^d$. For $d = 2$ the $n$- staircases are the Young diagrams of $n$. For a staircase $\lambda$, let $\bar{\lambda} := \mathbb{N}^d \setminus \lambda$ be its complement in $\mathbb{N}^d$ and let $\min(\bar{\lambda})$ be the unique finite set of coordinatewise minimal vectors in $\bar{\lambda}$. The $n$-staircases in $\mathbb{N}^d$ are in bijection with the monomial ideals in $\text{Hilb}_n^d$ via $I_\lambda := \text{ideal}\{x^v : v \in \min(\bar{\lambda})\}$.

Now fix any ideal $I \in \text{Hilb}_n^d$. An $n$-subset $\lambda \subset \mathbb{N}^d$ is *basic for* $I$ if the congruence classes modulo $I$ of the monomials $x^v$ with $v \in \lambda$ form a vector space basis for the quotient space $\mathbb{F}[x]/I$, or equivalently, if the $\mathbb{F}$-vector space $\lin\{x^v : v \in \lambda\}$ satisfies $\lin\{x^v : v \in \lambda\} \bigcap I = \{0\}$. If $\lambda$ is basic then the class $[f] = f + I$ of any $f \in \mathbb{F}[x]$ contains a unique representative in $\lin\{x^v : v \in \lambda\}$; let $[f]_\lambda$ denote this unique polynomial satisfying $[f]_\lambda \in \lin\{x^v : v \in \lambda\}$ and $f - [f]_\lambda \in I$.

A staircase $\lambda \in \binom{\mathbb{N}^d}{n}_{stair}$ is *initial for* $I$ if its monomial ideal $I_\lambda$ is the initial ideal $in_\prec(I) := \text{ideal}\{in_\prec(f) : f \in I\}$ of $I$ under some monomial order $\prec$. If $\lambda$ is initial then it is also basic and the unique reduced Gröbner basis of $I$ under $\prec$ is the set $G_\lambda(I) := \{x^u - [x^u]_\lambda : u \in \min(\bar{\lambda})\}$ consisting of precisely $|\min(\bar{\lambda})|$ polynomials. Let $\Lambda(I)$ denote the set of initial staircases of $I$. We shall need the following two propositions on basic sets and initial staircases of an ideal.

**Proposition 2.1** *Let $\prec$ be any monomial order and let $\lambda \in \Lambda(I)$ be the initial staircase of $I$ satisfying $I_\lambda = in_\prec(I)$. Then for any vector $u \in \mathbb{N}^d \setminus \lambda$ we have $[x^u]_\lambda \in lin\{x^v : v \in \lambda, \ v \prec u\}$.*

*Proof.* Let $f := x^u - [x^u]_\lambda \in I$ and let $x^v = in_\prec(f)$. Then $x^v \in in_\prec(I) = I_\lambda$ hence $v \notin \lambda$. Since $x^u$ is the only monomial in $f$ with exponent not in $\lambda$, it must be that $v = u$ hence $in_\prec(f) = x^u$, so all monomials involved in $[x^u]_\lambda$ are smaller than $x^u$ under $\prec$ as claimed. $\square$

**Proposition 2.2** *Let $\prec$ be any monomial order and let $\lambda \in \Lambda(I)$ be the initial staircase of $I \in \text{Hilb}_n^d$ satisfying $I_\lambda = in_\prec(I)$. Let $\mu$ be any basic set of $I$. Write $\lambda = \{\lambda_1, \ldots, \lambda_n\}$ and $\mu = \{\mu_1, \ldots, \mu_n\}$ with $\lambda_1 \prec \cdots \prec \lambda_n$ and $\mu_1 \prec \cdots \prec \mu_n$. Then $\lambda_k \preceq \mu_k$ for each $k = 1, \ldots, n$.*

*Proof.* Suppose indirectly $\mu_k \prec \lambda_k$ for some $k$, and let $U := \lin\{[x^{\mu_1}]_\lambda, \ldots, [x^{\mu_k}]_\lambda\}$. Then $\mu_i \preceq \mu_k \prec \lambda_k$ for all $i \leq k$ hence $[x^{\mu_i}]_\lambda \in \lin\{x^{\lambda_j} : j < k\}$ (clearly if $\mu_i \in \lambda$ and by Proposition 2.1 if $\mu_i \notin \lambda$). So $U \subseteq \lin\{x^{\lambda_j} : j < k\}$ hence $\dim(U) < k$. But $\mu$ is basic so $\{[x^{\mu_1}]_\lambda, \ldots, [x^{\mu_n}]_\lambda\}$ is linearly independent hence so is $\{[x^{\mu_1}]_\lambda, \ldots, [x^{\mu_k}]_\lambda\}$. Thus, $\dim(U) = k$, a contradiction. $\square$

Each vector $w$ in the nonnegative orthant $\mathbb{R}_+^d$ partially orders monomials $x^v$ by the value $w \cdot v$. For any generic $w$ this is a total order and hence a monomial order, and every monomial order arises that way from some generic $w \in \mathbb{R}_+^d$. The initial ideal of $I$ under $w$ is $in_w(I) := \text{ideal}\{in_w(f) : f \in I\}$ and is a monomial ideal if $w$ is generic. Declare two non-negative vectors $w$ and $w'$ equivalent if $in_w(I) = in_{w'}(I)$. The equivalence classes are relatively open convex cones forming a subdivision of $\mathbb{R}_+^d$ called the *Gröbner fan* of $I$ (cf. [8]). A vector $w$ lies in a full dimensional cone of the Gröbner fan if and only if $in_w(I)$ is a monomial ideal.

The (minimizing) *normal cone* of a face at a polyhedron $P$ in $\mathbb{R}^d$ is the relatively open cone of those vectors $w \in \mathbb{R}^d$ uniquely minimized over $P$ at that face. The collection of normal cones of all faces of $P$ is called the *normal fan* of $P$. Let $B$ be any polytope in $\mathbb{R}^d$ and let $P := B + \mathbb{R}_+^d$. Then the normal fan of the polyhedron $P$ forms a subdivision of $\mathbb{R}_+^d$ and the vertices of $P$ are precisely those vertices of $B$ whose normal cone at $B$ contains a strictly positive vector $w \in \mathbb{R}_+^d$.

A polyhedron $P$ in $\mathbb{R}^d$ is the *state polyhedron of the ideal $I$* (cf. [3]) if the Gröbner fan of $I$ equals the normal fan of $P$. This holds if and only if the set $\Lambda(I)$ of initial staircases is in bijection with the vertex set of $P$, with $w \in \mathbb{R}^d_+$ uniquely minimized over $P$ at its vertex corresponding to $\lambda \in \Lambda(I)$ if and only if the initial monomial ideal of $I$ under $w$ satisfies $in_w(I) = I_\lambda$.

We now describe a construction of the state polyhedron for which this bijection is very natural, in that the vertex corresponding to $\lambda \in \Lambda(I)$ is simply its vector sum $\sum \lambda \in \mathbb{N}^d$. Let $V_n^d := \bigcup \binom{\mathbb{N}^d}{n}_{stair}$ denote the union of all $n$-staircases in $\mathbb{N}^d$. Given an ideal $I \in \mathrm{Hilb}_n^d$, let $\Gamma(I) := \{\lambda \subset V_n^d : \lambda \text{ basic for } I\}$ denote the finite set of all $n$-subsets of $V_n^d$ basic for $I$. Since every initial staircase of $I$ is basic and is contained in $V_n^d = \bigcup \binom{\mathbb{N}^d}{n}_{stair}$ we have $\Lambda(I) \subseteq \Gamma(I)$. The following polytope will later enable the efficient computation of the state polyhedron.

**Definition 2.3** The *basis polytope* of $I \in \mathrm{Hilb}_n^d$ is the convex hull of sums of basic sets of $I$ in $V_n^d$,

$$\mathcal{B}(I) \quad := \quad \mathrm{conv}\{\sum \lambda \, : \, \lambda \in \Gamma(I)\,\} \quad \subset \quad \mathbb{R}^d \quad .$$

The state polyhedron of any $n$-long $d$-variate ideal is provided by the following theorem.

**Theorem 2.4** *The state polyhedron of $I \in \mathrm{Hilb}_n^d$ is provided by $\mathcal{S}(I) := B(I) + \mathbb{R}^d_+$. Furthermore:*

*(1) its vertex set is $\{\sum \lambda : \lambda \in \Lambda(I)\}$ and is in bijection with $\Lambda(I)$ via the map $\lambda \longrightarrow \sum \lambda$ ;*

*(2) a generic $w \in \mathbb{R}^d_+$ is minimized over $\mathcal{S}(I)$ and $\mathcal{B}(I)$ at $\sum \lambda$ with $\lambda \in \Lambda(I)$ and $I_\lambda = in_w(I)$ .*

*Proof.* Let $w \in \mathbb{R}^d_+$ be any generic vector and let $\lambda \in \Lambda(I)$ be the initial staircase of $I$ satisfying $in_w(I) = I_\lambda$. Consider any $\mu \in \Gamma(I)$. Writing $\lambda = \{\lambda_1, \ldots, \lambda_n\}$ and $\mu = \{\mu_1, \ldots, \mu_n\}$ with $w \cdot \lambda_1 < \cdots < w \cdot \lambda_n$ and $w \cdot \mu_1 < \cdots < w \cdot \mu_n$, we find by Proposition 2.2 that for each $k$ we have $w \cdot \lambda_k \leq w \cdot \mu_k$, with equality if and only if $\lambda_k = \mu_k$ since $w$ yields a total order. Thus,

$$w \cdot \sum \lambda \quad = \quad \sum_{k=1}^n w \cdot \lambda_k \quad \leq \quad \sum_{k=1}^n w \cdot \mu_k \quad = \quad w \cdot \sum \mu$$

with equality if and only if $\lambda_k = \mu_k$ for all $k$. We find that for all $\mu \in \Gamma(I)$ other than $\lambda$ we have $w \cdot \sum \lambda < w \cdot \sum \mu$ hence $w$ is uniquely minimized over $\mathcal{B}(I) + \mathbb{R}^d_+$ at $\sum \lambda$. This shows in particular that $\sum \lambda$ is indeed a vertex of $\mathcal{S}(I)$ and $\mathcal{B}(I)$ and that $\sum \mu \neq \sum \lambda$ for any $\mu \in \Gamma(I)$ other than $\lambda$. Since any $\lambda \in \Lambda(I)$ satisfies $I_\lambda = in_w(I)$ for some generic $w \in \mathbb{R}^d_+$, the map $\lambda \to \sum \lambda$ is indeed a bijection from $\Lambda(I)$ onto $\{\sum \lambda : \lambda \in \Lambda(I)\}$. Since generic vectors are dense in $\mathbb{R}^d_+$ and a fan is determined by its maximal dimensional cones we find that the Gröbner fan of $I$ equals the normal fan of $\mathcal{S}(I) = \mathcal{B}(I) + \mathbb{R}^d_+$, showing that $\mathcal{S}(I)$ is indeed the state polyhedron of $I$. □

Theorem 2.4 remains valid if we replace $\mathcal{B}(I)$ by the convex hull $\mathcal{B} = \mathrm{conv}\{\sum \lambda \, : \, \lambda \in \Gamma\}$ with $\Gamma \supseteq \Lambda(I)$ any collection of basic sets of $I$ which contains all initial staircases of $I$. In particular, it holds with $\Gamma = \Lambda(I)$ and with $\Gamma = \{\lambda \in \binom{\mathbb{N}^d}{n}_{stair} : \lambda \text{ basic for } I\}$. However, these choices do not lend themselves to efficient algorithmic construction of the state polyhedron and the universal Gröbner basis. As we shall see, we need to take $\Gamma$ to be the collection of all basic sets of $I$ contained in some subset $V \subseteq \mathbb{N}^d$. The smallest such set yielding $\Gamma$ which contains all $n$-staircases in $\mathbb{N}^d$ is their union $V_n^d$, leading to our choice of $\Gamma(I)$ and $\mathcal{B}(I)$.

Let $\mathcal{U}(I)$ denote the *universal Gröbner basis of $I$* defined as the union of all reduced Gröbner bases of $I$ and hence simultaneously providing a Gröbner basis for $I$ under any monomial order. Theorem 2.4 yields the following polynomial upper bounds on the number of distinct reduced Gröbner bases and the size of the universal Gröbner basis of any ideal in $\mathrm{Hilb}_n^d$. The proof is similar to that of [10] for the case of vanishing ideals of point configurations.

**Corollary 2.5** *For every fixed $d$, the following hold for any $n$-long $d$-variate ideal $I \in \mathrm{Hilb}_n^d$:*

*(1) the number of distinct reduced Gröbner bases of $I$ is $O(n^{2d\frac{d-1}{d+1}})$;*

*(2) the number of elements in the universal Gröbner basis $\mathcal{U}(I)$ is $O(n^{2d-3+\frac{3d-1}{d(d+1)}})$.*

*Proof.* Consider any $\lambda \in \binom{\mathbb{N}^d}{n}_{stair}$. Then, for $i = 1, \ldots, d$, each $v \in \lambda$ satisfies $v_i < n$ hence the $i$th coordinate value of $\sum \lambda$ is less than $n^2$. Thus, the lattice polytope $P = \mathrm{conv}\{\sum \lambda : \lambda \in \Lambda(I)\}$ is contained in the cube $[0, n^2]^d$ and hence has $O(\mathrm{vol}(P)^{\frac{d-1}{d+1}}) = O(n^{2d\frac{d-1}{d+1}})$ vertices (cf. [1]). By Theorem 2.4, the state polyhedron $\mathcal{S}(I) = P + \mathbb{R}_+^d$ and $P$ have the same vertex set, giving (1).

Next, for any $\lambda \in \binom{\mathbb{N}^d}{n}_{stair}$, the size of the set $\min(\bar{\lambda})$ of minimal elements not in $\lambda$ is $O(n^{\frac{d-1}{d}})$ (cf. [4]). Since the reduced Gröbner basis $G_\lambda$ of $I$ corresponding to $\lambda$ has $|\min(\bar{\lambda})|$ elements, the product of the bound on any $|\min(\bar{\lambda})|$ and the bound just established in (1) on the number of reduced Gröbner basis yields the bound in (2) on the size of their union $\mathcal{U}(I)$. □

## 3 The Hilbert zonotope and its universality

While the number of initial staircases of $I \in \mathrm{Hilb}_n^d$ is polynomial in $n$ for any fixed $d$ by Corollary 2.5, the cardinality of $\Gamma(I)$ is typically exponential in $n$ even for $d = 2$; for instance, if $I$ is the bivariate vanishing ideal of any generic $n$ points in the plane over an infinite field (see discussion of point configurations in Section 5), then all $n$-subsets of $V_n^d$ are basic, and hence so are all $n$-staircases in $\mathbb{N}^2$ which are in bijection with number partitions of $n$. Thus, it is not possible to filter $\Lambda(I)$ out of $\Gamma(I)$ and construct $\mathcal{B}(I)$ or $\mathcal{S}(I)$ directly in polynomial time.

To overcome this we now introduce, for each pair of positive integers $d$ and $n$, the *Hilbert zonotope* $\mathcal{H}_n^d$. As we shall see, it is *universal* for the Hilbert scheme $\mathrm{Hilb}_n^d$ in that it provides a refinement of the basis polytope $\mathcal{B}(I)$ and the state polyhedron $\mathcal{S}(I)$ of every ideal $I \in \mathrm{Hilb}_n^d$.

Recall that $V_n^d = \bigcup \binom{\mathbb{N}^d}{n}_{stair}$ is the union of all $n$-staircases in $\mathbb{N}^d$. Call an element of the symmetrization $V_n^d - V_n^d = \{u - v : u, v \in V_n^d\}$ *primitive* if it is not a nonnegative integer multiple of another element of $V_n^d - V_n^d$. For $n \geq 2$ let $D_n^d$ be the set of primitive elements of $V_n^d - V_n^d$, and for $n = 1$ let $D_1^d := \pm\{e_1, \ldots, e_d\}$. We make the following fundamental definition.

**Definition 3.1** The *Hilbert zonotope* $\mathcal{H}_n^d$ is the following Minkowski sum of line segments,

$$\mathcal{H}_n^d \quad := \quad \sum_{v \in D_n^d} [0,1] \cdot v \quad \subset \quad \mathbb{R}^d \quad .$$

Note that $D_n^d = -D_n^d$ is centrally symmetric, hence so is $\mathcal{H}_n^d$ which equivalently can be defined as $\sum [-1,1] \cdot v$ by summing over only one of each pair $\{-v, v\}$ of antipodal primitive elements.

The Hilbert zonotope is well behaved in the following sense.

**Proposition 3.2** *Fix any $d$. Then the number of vertices of the Hilbert zonotope $\mathcal{H}_n^d$ satisfies $O(n^{2(d-1)} (\log n)^{2(d-1)^2})$; further, in polynomial time using that many arithmetic operations, all its vertices can be listed, each $h$ along with a vector $w(h)$ uniquely minimized over $\mathcal{H}_n^d$ at $h$.*

*Proof.* First note that $V_n^d = \bigcup \binom{\mathbb{N}^d}{n}_{stair}$ is given by $V_n^d = \{v \in \mathbb{N}^d : \prod_{i=1}^d (v_i + 1) \leq n\}$: indeed, if $v \in \mathbb{N}^d$ lies in some $n$-staircase $\lambda$ then the entire box $\{u \in \mathbb{N}^d : u \leq v\}$ is contained in $\lambda$ hence $\prod_{i=1}^d (v_i + 1) = |\{u \in \mathbb{N}^d : u \leq v\}| \leq |\lambda| = n$; conversely, if $v \in \mathbb{N}^d$ satisfies $\prod_{i=1}^d (v_i + 1) \leq n$ then the box above can be augmented with multiples of unit vectors to an $n$-staircase $\lambda$ containing $v$. Now the cardinality of $\{v \in \mathbb{N}^d : \prod_{i=1}^d (v_i + 1) \leq n\}$ obeys the bound $O(n (\log n)^{d-1})$ (cf. [12]). Thus, the Hilbert zonotope $\mathcal{H}_n^d$ is the Minkowski sum of $N := \frac{1}{2}|D_n^d| = O(|V_n^d|^2) = O(n^2 (\log n)^{2(d-1)})$ line segments and therefore (cf. [6, 9]) has $O(N^{d-1})$ vertices which can all be enumerated, each $h$ along with a vector $w(h)$ uniquely minimized at $h$, using $O(N^{d-1})$ arithmetic operations, giving the claimed bounds. □

**Example 3.3** We compute the Hilbert zonotope $\mathcal{H}_3^2$ for $d = 2$, $n = 3$. Using Proposition 3.2 we find the set $H$ of vertices of $\mathcal{H}_3^2$ which shows that it is a (centrally symmetric) 10-gon, as well as the corresponding set $W$ of the linear functionals uniquely minimized at its various vertices:

$$V_3^2 = \left\{ \begin{bmatrix} 0 \\ 0 \end{bmatrix}, \begin{bmatrix} 1 \\ 0 \end{bmatrix}, \begin{bmatrix} 0 \\ 1 \end{bmatrix}, \begin{bmatrix} 2 \\ 0 \end{bmatrix}, \begin{bmatrix} 0 \\ 2 \end{bmatrix} \right\}, \quad D_3^2 = \pm \left\{ \begin{bmatrix} 1 \\ 0 \end{bmatrix}, \begin{bmatrix} 0 \\ 1 \end{bmatrix}, \begin{bmatrix} 1 \\ -1 \end{bmatrix}, \begin{bmatrix} 1 \\ -2 \end{bmatrix}, \begin{bmatrix} 2 \\ -1 \end{bmatrix} \right\},$$

$$H = \pm \left\{ \begin{bmatrix} -5 \\ 5 \end{bmatrix}, \begin{bmatrix} -5 \\ 3 \end{bmatrix}, \begin{bmatrix} -3 \\ -1 \end{bmatrix}, \begin{bmatrix} -1 \\ -3 \end{bmatrix}, \begin{bmatrix} 3 \\ -5 \end{bmatrix} \right\}, W = \pm \left\{ \begin{bmatrix} 1 \\ -1 \end{bmatrix}, \begin{bmatrix} 3 \\ 1 \end{bmatrix}, \begin{bmatrix} 3 \\ 2 \end{bmatrix}, \begin{bmatrix} 2 \\ 3 \end{bmatrix}, \begin{bmatrix} 1 \\ 3 \end{bmatrix} \right\}.$$

We need to recall a few facts about matroids and matroid polytopes. Let $M = (V, \Gamma)$ be a matroid over a finite set $V$ with collection of bases $\Gamma \subseteq 2^V$. Its *matroid polytope* is defined as the convex hull $\mathcal{B}(M) := \text{conv}\{\mathbf{1}_B : B \in \Gamma\} \subset \mathbb{R}^V$, where $\mathbf{1}_B$ denotes the *incidence vector* of $B \subseteq V$, that is, the $\{0,1\}$ vector in $\mathbb{R}^V$ with support $B$. This is a well known object of importance in combinatorial optimization. Below, let $e_v \in \mathbb{R}^V$ denote the unit vector indexed by $v \in V$.

**Proposition 3.4** *Every 1-face of the matroid polytope is equal to $e_u - e_v$ for some $u, v \in V$.*

*Proof.* Consider any pair $A, B \in \Gamma$ of bases such that $[\mathbf{1}_A, \mathbf{1}_B]$ is an edge (that is, a 1-face) of $\mathcal{B}(M)$, and let $w \in \mathbb{R}^V$ be a linear functional uniquely maximized over $\mathcal{B}(M)$ at that edge. If $A \setminus B = \{u\}$ is a singleton then $B \setminus A = \{v\}$ is a singleton as well in which case $\mathbf{1}_A - \mathbf{1}_B = e_u - e_v$ and we are done. Suppose then, indirectly, that it is not, and pick an element $u$ in the symmetric difference $A \Delta B := (A \setminus B) \cup (B \setminus A)$ of $A$ and $B$ of minimum value $w_u$. Without loss of generality assume $u \in A \setminus B$. Then there is a $v \in B \setminus A$ such that $C := A \setminus \{u\} \cup \{v\}$ is a basis of $M$. Since $|A \Delta B| > 2$, $C$ is neither $A$ nor $B$. By the choice of $u$, this basis satisfies $w \cdot \mathbf{1}_C = w \cdot \mathbf{1}_A - w_u + w_v \geq w \cdot \mathbf{1}_A$, hence also a maximizer of $w$ over $\mathcal{B}(M)$, a contradiction. □

A polyhedron $P$ is a *refinement* of a polyhedron $Q$ if the normal fan of $P$ is a refinement of that of $Q$, that is, the closure of each normal cone of $Q$ is the union of closures of normal cones of $P$. The significance of the Hilbert zonotope is now demonstrated by the following theorem.

**Theorem 3.5** *The Hilbert zonotope $\mathcal{H}_n^d$ is a refinement of both the basis polytope $\mathcal{B}(I)$ and the state polyhedron $\mathcal{S}(I)$ of every member $I$ of the Hilbert scheme $\mathrm{Hilb}_n^d$ of n-long d-variate ideals.*

*Proof.* Consider any $I \in \mathrm{Hilb}_n^d$. Let $M := (V_n^d, \Gamma(I))$ be the matroid over $V_n^d$ with collection of bases $\Gamma(I)$, which is the restriction to $V_n^d$ of the infinite matroid over $\mathbb{N}^d$ of all basic sets of $I$. Let $\mathcal{B}(M) := \mathrm{conv}\{\mathbf{1}_\lambda : \lambda \in \Gamma(I)\} \subset \mathbb{R}^{V_n^d}$ be the matroid polytope of $M$ and let

$$\pi : \mathbb{R}^{V_n^d} \longrightarrow \mathbb{R}^d \ : \ e_v \mapsto v$$

be the natural projection sending the unit vector $e_v$ corresponding to $v \in V_n^d$ to the vector $v \in \mathbb{N}^d \subset \mathbb{R}^d$. Then for each $\lambda \in \Gamma(I)$ we have $\pi(\mathbf{1}_\lambda) = \sum \lambda$ hence the basis polytope of $I$

$$\mathcal{B}(I) \ = \ \mathrm{conv}\{\sum \lambda : \lambda \in \Gamma(I)\} \ = \ \mathrm{conv}\{\pi(\mathbf{1}_\lambda) : \lambda \in \Gamma(I)\} \ = \ \pi(\mathcal{B}(M))$$

is a projection of the matroid polytope. Thus, each edge of $\mathcal{B}(I)$ is the projection of some edge of $\mathcal{B}(M)$ hence, by Proposition 3.4, is equal to $\pi(e_u - e_v) = u - v$ for some pair $u, v \in V_n^d$ and therefore parallel to some element in $D_n^d$. Thus, the Hilbert zonotope $\mathcal{H}_n^d = \sum_{v \in D_n^d} [0, 1] \cdot v$ is the Minkowski sum of a set of segments containing all edge directions of $\mathcal{B}(I)$ and therefore its normal fan is a refinement of the normal fan of $\mathcal{B}(I)$.

Next, consider any face $F$ of $\mathcal{S}(I)$. Then $F$ is also a face of $\mathcal{B}(I) \subset \mathcal{S}(I)$ and hence the closure of its normal cone at $\mathcal{B}(I)$ is a union of closures of normal cones of $\mathcal{H}_n^d$. But $D_n^d$ contains all unit vectors which implies that each normal cone of $\mathcal{H}_n^d$ is contained in the interior of some orthant. Thus, the closure of the normal cone of $F$ at $\mathcal{S}(I)$ is the union of the closures of those normal cones of $\mathcal{H}_n^d$ contained in the normal cone of $F$ at $\mathcal{B}(I)$ which lie in the nonnegative orthant $\mathbb{R}_+^d$. □

Let us call the coarsest common refinement of the state polyhedra of all ideals $I$ on the Hilbert scheme the *Hilbert polytope*. As pointed out to us by D. Bayer, the Hilbert zonotope may in general be finer than the Hilbert polytope; in particular, it may have more vertices. However, the Hilbert zonotope allows efficient algorithmic treatment while the Hilbert polytope might not.

## 4  Computing the state polyhedron and universal Gröbner basis

We now use the Hilbert zonotope to efficiently compute the set of initial staircases, the state polyhedron and the universal Gröbner basis of any ideal on the Hilbert scheme.

Let $U_n^d := \{u + e_i : u \in V_n^d,\ 0 \leq i \leq d\}$ with $e_0 := 0$ and $e_i$ the $i$th unit vector for $1 \leq i \leq d$. Then $U_n^d$ contains, along with every $n$-staircase $\lambda$, all vectors $u + e_i$ with $u$ any vector in $\lambda$ and $e_i$ any unit vector, hence also $\min(\bar{\lambda})$. On the other hand, $U_n^d \subseteq V_{2n}^d$ and hence, for fixed $d$, obeys the same upper bound $O(n(\log n)^{d-1})$ on its cardinality as does $V_{2n}^d$.

We assume that $I \in \mathrm{Hilb}_n^d$ is presented by its Gröbner basis $G_\lambda := \{x^u - [x^u]_\lambda : u \in \min(\bar{\lambda})\}$ under some monomial order $\prec$. Such a presentation, say with respect to the degree reverse lexicographic order, is known to be efficiently computable from any generating set - see discussion at the end of this section. Given such an ideal we will need the representative $[x^u]_\lambda$ for every $u$ in the set $U_n^d$. We include the short proof of the following adaptation of [5, Propoistion 3.1].

**Proposition 4.1** *Fix any $d$. Then, given any ideal $I = \mathrm{ideal}(G_\lambda) \in \mathrm{Hilb}_n^d$, the representatives $[x^u]_\lambda = \sum_{v \in \lambda} a_{v,u} \cdot x^v$ of all $u \in U_n^d$ can be computed using $O(n^3 (\log n)^{d-1})$ arithmetic operations.*

*Proof.* Compute $[x^u]_\lambda$ for the elements $u \in U_n^d$ in the order $\prec$ as follows. If $u \in \lambda$ then $[x^u]_\lambda = x^u$. If $u \in \min(\bar{\lambda})$ then $x^u - [x^u]_\lambda$ is in the given reduced Gröbner basis $G_\lambda$ from which $[x^u]_\lambda$ can be recovered. Otherwise $s := u - e_i \in U_n^d \setminus \lambda$ for some $i$. By Proposition 2.1 $[x^s]_\lambda \in \mathrm{lin}\{x^t : t \in T\}$ with $T := \{t \in \lambda : t \prec s\}$. Now $s \prec u$ and $t + e_i \prec s + e_i = u$ for all $t \in T$, hence $[x^s]_\lambda = \sum_{t \in T} a_{t,s} x^t$ and $[x^{t+e_i}]_\lambda = \sum_{v \in \lambda} a_{v,t+e_i} x^v$ for all $t \in T$ are already available. Thus

$$[x^u]_\lambda = [x^{s+e_i}]_\lambda = \sum_{t \in T} a_{t,s}[x^{t+e_i}]_\lambda = \sum_{t \in T} a_{t,s} \sum_{v \in \lambda} a_{v,t+e_i} x^v = \sum_{v \in \lambda} (\sum_{t \in T} a_{t,s} a_{v,t+e_i}) x^v$$

is now obtained using $O(n^2)$ arithmetic operations. Since $|U_n^d| = O(n(\log n)^{d-1})$ we are done.  □

We are finally in position to provide the efficient procedure for computing the set of initial staircases, the state polyhedron and the universal Gröbner basis of any ideal on the Hilbert scheme. Theorem 3.5 enables us to bypass the difficulty caused by the exponential size of $\Gamma(I)$.

**Theorem 4.2** *For every fixed $d$ there is a polynomial time algorithm that, given any $n$-long $d$-variate ideal $I \in \mathrm{Hilb}_n^d$, computes its set $\Lambda(I)$ of initial staircases, its state polyhedron $\mathcal{S}(I)$, and its universal Gröbner basis $\mathcal{U}(I)$ using $O(n^{2d+1} (\log n)^{(2d-1)(d-1)})$ arithmetic operations.*

*Proof.* First enumerate, as in Proposition 3.2 and hence within the claimed complexity bound, the $O(n^{2(d-1)} (\log n)^{2(d-1)^2})$ vertices of $\mathcal{H}_n^d$, each vertex $h$ along with a vector $w(h)$ uniquely minimized over $\mathcal{H}_n^d$ at $h$. We claim that any (coordinatewise) positive $w(h)$ on the list is uniquely minimized over $\mathcal{B}(I)$ at $\sum \mu$ for some initial staircase $\mu \in \Lambda(I)$ and, conversely, for every $\mu \in \Lambda(I)$ there is some (possibly many) positive $w(h)$ on the list uniquely minimized over $\mathcal{B}(I)$ at $\sum \mu$. Consider first any positive $w(h)$ on the list. Since $\mathcal{H}_n^d$ refines $\mathcal{B}(I)$ by Theorem 3.5, the vector $w(h)$ and hence a positive generic perturbation of it lie in the normal cone of some vertex of $\mathcal{B}(I)$.

By Theorem 2.4, this vertex is $\sum \mu$ for some initial staircase $\mu \in \Lambda(I)$. Conversely, consider any initial staircase $\mu \in \Lambda(I)$. Then, by Theorem 2.4, $\sum \mu$ is a vertex of $\mathcal{S}(I)$ hence its normal cone is contained in the interior of the nonnegative orthant. Since $\mathcal{H}_n^d$ refines $\mathcal{S}(I)$ by Theorem 3.5, this normal cone contains the normal cone of some (possibly many) vertex $h$ of $\mathcal{H}_n^d$. The vector $w(h)$ of that vertex $h$ on the list must then be positive.

We proceed to show that for each positive $w(h)$ on the list we can efficiently compute the minimizing vertex $\sum \mu$ of $\mathcal{B}(I)$, the initial staircase $\mu$, and the reduced Gröbner basis $G_\mu$.

First, compute, as in Proposition 4.1, using the given presentation $I = \mathrm{ideal}(G_\lambda)$ of $I$ by some reduced Gröbner basis, the representative $[x^u]_\lambda = \sum_{v \in \lambda} a_{v,u} \cdot x^v$ of every vector $u \in U_n^d$.

Now pick any vertex $h$ of $\mathcal{H}_n^d$ with $w := w(h)$ positive. Use the following greedy algorithm to find a basic set $\mu = \{\mu_1, \ldots, \mu_n\} \in \Gamma(I)$: for $i = 1, \ldots, n$ pick $\mu_i$ to be an element in $V_n^d \subset U_n^d$ of minimal value $w \cdot \mu_i$ with the property that $\{[x^{\mu_j}]_\lambda : j \leq i\}$ is linearly independent. Then for any other basic set $\nu \in \Gamma(I)$, writing $\nu = \{\nu_1, \ldots, \nu_n\}$ with $w \cdot \nu_1 < \cdots < w \cdot \nu_n$, we have $w \cdot \mu_k \leq w \cdot \nu_k$ for all $k$ hence $w \cdot \sum \mu \leq w \cdot \sum \nu$. This shows that $\sum \mu$ minimizes $w$ over $\mathcal{B}(I)$ hence, as shown in the first paragraph of this proof, is the unique minimizing vertex of the positive $w(h) = w$ over $\mathcal{B}(I)$, and $\mu \in \Lambda(I)$ is the corresponding initial staircase.

The determination of $\mu$ by the greedy algorithm can be efficiently implemented by performing Gaussian elimination on the fly as follows. Totally order $U_n^d$ compatibly with $w$ so that if $u$ precedes $v$ then $w \cdot u \leq w \cdot v$ and let $\lambda = \{\lambda_1, \ldots, \lambda_n\}$ be an arbitrary labeling of $\lambda$. Thus, the coefficients $a_{v,u}$ for all $u \in U_n^d$ and $v \in \lambda$ now form an $n \times |U_n^d|$ matrix $A$ over $\mathbb{F}$. For each $u \in U_n^d$ let $A^u \in \mathbb{F}^n$ denote the column of that matrix corresponding to $u$, with $A_i^u := a_{\lambda_i, u}$. Now, in the greedy algorithm, for $i = 1, \ldots, n$, pick $\mu_i$ to be the first element in $V_n^d \subset U_n^d$ whose column $A^{\mu_i}$ contains a nonzero coordinate $A_k^{\mu_i}$ for some $k \geq i$. Apply suitable row operations to $A$ so as to transform $A^{\mu_i}$ to the unit vector $e_i$ while maintaining $A^{\mu_j} = e_j$ for all $j < i$. This consumes $O(n \cdot |U_n^d|)$ arithmetic operations per $\mu_i$, totaling to $O(n^2 \cdot |U_n^d|)$ operations.

Moreover, the updated matrix $A$ at the end of this iterated process is the matrix of coordinates in the new basis $\mu$: for each $u \in U_n^d$, the representative $[x^u]_\mu$ is simply read off from the corresponding column as $[x^u]_\mu = \sum_{i=1}^n A_i^u x^{\mu_i}$. Now, the set $\min(\bar{\mu})$ consists precisely of those $u \in U_n^d \setminus \mu$ with the property that for each $i = 1, \ldots, d$ either $u_i = 0$ or $u - e_i \in \mu$, hence can be quickly filtered out of $U_n^d$. So, within the same complexity bound of $O(n^2 \cdot |U_n^d|)$ operations we obtain the reduced Gröbner basis $G_\mu := \{x^u - [x^u]_\mu : u \in \min(\bar{\mu})\}$ corresponding to $\mu$.

Summarizing, for each positive $w(h)$ on the list of vertices $h$ of $\mathcal{H}_n^d$, we can find the initial staircase $\mu \in \Lambda(I)$ for which $\sum \mu$ is the unique minimizer of $w(h)$ over $\mathcal{B}(I)$ and the corresponding reduced Gröbner basis $G_\mu$ using $O(n^2 \cdot |U_n^d|)$ operations. Therefore, the entire set $\Lambda(I)$ of initial staircases, the state polyhedron $\mathcal{S}(I) = \mathrm{conv}\{\sum \mu : \mu \in \Lambda(I)\} + \mathbb{R}_+^d$, and the universal Gröbner basis $\mathcal{U}(I) = \bigcup_{\mu \in \Lambda(I)} G_\mu$ of the given ideal $I$ can be produced using $O(n^{2(d-1)} (\log n)^{2(d-1)^2} \cdot n^2 \cdot n (\log n)^{d-1})$ arithmetic operations as claimed. □

**Example 4.3** We compute the universal Gröbner basis of the ideal $I = \mathrm{ideal}(G_\lambda) \in \mathrm{Hilb}_3^2$ for $d = 2$, $n = 3$ over $\mathbb{F} = \mathbb{R}$, presented by $G_\lambda = \{x_1^3 - 3x_1^2 + 3x_1 - 1, x_2 - x_1 + 1\}$ which is its reduced Gröbner basis under the lexicographic order with $x_2 > x_1$, with $\lambda$ the staircase $\lambda = \{00, 10, 20\}$. First, we obtain the set of $w \in \mathbb{R}_+^2$ via the Hilbert zonotope $\mathcal{H}_3^2$ as in Example 3.3, which is $W_+ := \{31, 32, 23, 13\}$. Next, we compute the $3 \times 10$ matrix $A$ of coefficients of the representatives $[x^u]_\lambda = a_{1,u} \cdot 1 + a_{2,u} \cdot x_1 + a_{3,u} \cdot x_1^2$ for all $u \in U_3^2$, whose rows and columns are indexed by $\lambda$ and $U_3^2$ respectively, and obtain

$$A = \begin{array}{c} \\ 00 \\ 10 \\ 20 \end{array} \begin{pmatrix} 00 & 10 & 20 & 30 & 01 & 11 & 21 & 02 & 12 & 03 \\ 1 & 0 & 0 & 1 & -1 & 0 & 1 & 1 & 1 & 0 \\ 0 & 1 & 0 & -3 & 1 & -1 & -3 & -2 & -2 & 0 \\ 0 & 0 & 1 & 3 & 0 & 1 & 2 & 1 & 1 & 0 \end{pmatrix} .$$

Now consider some vector in $W_+$, say $w := 32$. Reorder $U_3^2$ compatibly with $w$ and suitably permute the columns of $A$. Apply the greedy algorithm and find the new initial staircase $\mu = \{00, 01, 02\}$. Next apply suitable row operations to $A$ to make the transformation to the new basis and obtain the following updated matrix, with rows and columns suitably re-labeled,

$$A = \begin{array}{c} \\ 00 \\ 01 \\ 02 \end{array} \begin{pmatrix} 00 & 01 & 10 & 02 & 11 & 03 & 20 & 12 & 21 & 30 \\ 1 & 0 & 1 & 0 & 0 & 0 & 0 & 0 & 0 & 1 \\ 0 & 1 & 1 & 0 & 1 & 0 & 2 & 0 & 1 & 3 \\ 0 & 0 & 0 & 1 & 1 & 0 & 1 & 1 & 2 & 3 \end{pmatrix} .$$

Now, we have $\min(\bar{\mu}) = \{10, 03\}$ and so the new reduced Gröbner basis is read off from the third and sixth columns of the new matrix to be $G_\mu = \{x_1 - x_2 - 1, x_2^3\}$.

Repeating this for each of the other three vectors in $W_+$, we keep getting either $\lambda$ or $\mu$. We conclude that here $\Lambda(I) = \{\lambda, \mu\}$ consists of two staircases only, the state polyhedron is $\mathcal{S}(I) = \mathrm{conv}\{\sum \lambda, \sum \mu\} + \mathbb{R}_+^2 = [30, 03] + \mathbb{R}_+^2$, and the universal Gröbner basis is

$$\mathcal{U}(I) \quad = \quad G_\lambda \cup G_\mu \quad = \quad \{x_1^3 - 3x_1^2 + 3x_1 - 1,\ x_2 - x_1 + 1,\ x_1 - x_2 - 1,\ x_2^3\} \quad .$$

It is known (cf. [5] and references therein) that the reduced Gröbner basis under the degree reverse lexicographic order of any ideal $I = \mathrm{ideal}(F) \in \mathrm{Hilb}_n^d$ presented by any set of generators can be computed in time $O(D^{d^2})$ where $D$ is the maximal degree of any generator $f \in F$. Thus, Theorem 4.2 also implies, for any fixed $d$, an efficient algorithm for computing the state polyhedron and the universal Gröbner basis of any ideal of the Hilbert scheme $\mathrm{Hilb}_n^d$ presented by any set of generators using a number of arithmetic operations polynomial in $n$ and $D$.

We conclude by pointing out that Theorems 3.5 and 4.2 yield the practical outcome that, for each $d$ and $n$, a list $W_n^d$ of positive $w(h)$ of vertices $h$ of $\mathcal{H}_n^d$ can be computed once and for all, providing a *universal set of monomial orders* for the Hilbert scheme $\mathrm{Hilb}_n^d$ that allows the efficient computation of the universal Gröbner basis $\mathcal{U}(I)$ of any given $n$-long $d$-variate ideal $I$.

## 5   Examples: point configurations and lattice ideals

We now interpret some of the notions and results discussed above for several special classes of ideals of the Hilbert scheme $\mathrm{Hilb}_n^d$. We start with the simple class of monomial ideals.

**Monomial Ideals.** Recall that the monomial ideals in $\text{Hilb}_n^d$ are in bijection with the $n$-staircases in $\mathbb{N}^d$ via $I_\lambda = \text{ideal}\{x^u : u \in \min(\bar{\lambda})\}$. Consider any such ideal $I_\lambda$. Then $\Gamma(I_\lambda) = \{\lambda\}$ is a singleton: indeed, if $\mu$ is any other $n$-subset of $\mathbb{N}^d$ and $u \in \mu \setminus \lambda$ then $x^u \in \text{lin}\{x^v : v \in \mu\} \bigcap I_\lambda$ hence $\mu$ is not basic. For any monomial order $\prec$ we have $in_\prec(I_\lambda) = I_\lambda$ hence $\Lambda(I_\lambda) = \{\lambda\}$ as well. Thus, $\mathcal{B}(I_\lambda) = \{\sum \lambda\}$ is a single point and the state polyhedron is $\mathcal{S}(I_\lambda) = \{\sum \lambda\} + \mathbb{R}_+^d$. For every $u \in \mathbb{N}^d \setminus \lambda$ we have $x^u \in I_\lambda$ hence $[x^u]_\lambda = 0$. Thus, the universal Gröbner basis of $I$ equals the unique reduced Gröbner basis $G_\lambda$ and are both given by $\mathcal{U}(I_\lambda) = G_\lambda = \{x^u : u \in \min(\bar{\lambda})\}$.

**Point Configurations.** The vanishing ideal $I_C := \{f \in \mathbb{F}[x] : f(c_1) = \cdots = f(c_n) = 0\}$ of a configuration $C = \{c_1, \ldots, c_n\}$ of $n$ distinct points in affine space $\mathbb{F}^d$ is a radical ideal of length $\dim(\mathbb{F}[x]/I_C) \leq n$. Assume throughout this example that $\mathbb{F}$ is infinite, which implies that the length of $I_C$ is exactly $n$ and hence $I_C \in \text{Hilb}_n^d$. For $\lambda = \{\lambda_1, \ldots, \lambda_n\} \subset \mathbb{N}^d$ let

$$C^\lambda := \begin{pmatrix} c_1^{\lambda_1} & c_1^{\lambda_2} & \cdots & c_1^{\lambda_n} \\ c_2^{\lambda_1} & c_2^{\lambda_2} & \cdots & c_2^{\lambda_n} \\ \vdots & \vdots & \ddots & \vdots \\ c_n^{\lambda_1} & c_n^{\lambda_2} & \cdots & c_n^{\lambda_n} \end{pmatrix}, \quad \text{where} \quad c_i^{\lambda_j} = \prod_{k=1}^d c_{i,k}^{\lambda_{j,k}}.$$

Then $\lambda$ is basic for $I_C$ if and only if $C^\lambda$ is nonsingular: indeed, the vector $(f(c_1), \ldots, f(c_n)) \in \mathbb{F}^n$ of evaluations of a polynomial $f = \sum_{i=1}^n a_i x^{\lambda_i} \in \text{lin}\{x^{\lambda_1}, \ldots, x^{\lambda_n}\}$ at the points of $C$ is provided by $C^\lambda \cdot A$ and is the zero vector in $\mathbb{F}^n$ if and only if $f$ lies in $I$; thus, $\text{lin}\{x^{\lambda_1}, \ldots, x^{\lambda_n}\} \cap I_C = \{0\}$ if and only if $C^\lambda$ is nonsingular. Thus, $\Gamma(I_C) = \{\lambda \subset V_n^d : \det(C^\lambda) \neq 0\}$; the basis polytope is $\mathcal{B}(I_C) = \text{conv}(\sum \lambda : \lambda \in \Gamma(I_C))$; the state polyhedron is $\mathcal{S}(I_C) = \mathcal{B}(I_C) + \mathbb{R}_+^d$; and the set of initial staircases is $\Lambda(I_C) = \{\lambda \in \Gamma(I_C) : \sum \lambda \text{ vertex of } \mathcal{S}(I_C)\}$. For every sufficiently generic configuration, say one satisfying $\det(C^\lambda) \neq 0$ for all $\lambda \subset V_n^d$, the state polyhedron $\mathcal{S}(I_C)$ and the set of initial staircases $\Lambda(I_C)$ coincide, respectively, with the *corner cut* polyhedron $P_n^d$ and the set $\binom{\mathbb{N}^d}{n}_{cut}$ of $n$-element *corner cuts* in $\mathbb{N}^d$ introduced and studied in [10]. If $\prec$ is a monomial order and $\lambda \in \Lambda(I_C)$ is the initial staircase with $in_\prec(I_C) = I_\lambda$ then for every $u \in \mathbb{N}^d$ we have $[x^u]_\lambda = x^u - \det((\{x\} \cup C)^{(\{u\} \cup \lambda)}) \cdot \det^{-1}(C^\lambda)$. So the reduced Gröbner basis of $I_C$ under $\prec$ is

$$G_\lambda = \{\det((\{x\} \cup C)^{(\{u\} \cup \lambda)}) \cdot \det^{-1}(C^\lambda) : u \in \min(\bar{\lambda})\}$$

and the universal Gröbner basis $\mathcal{U}(I_C) = \cup_{\lambda \in \Lambda(I_C)} G_\lambda$ is efficiently computable by our methods. This extends the results of [10] where an efficient construction of $\mathcal{U}(I_C)$ for *generic* configurations only was provided, based on the separable-partitions methods of [2, 7].

**Lattice Ideals.** The binomial ideal of an integer lattice $L \subseteq \mathbb{Z}^d$ is $I_L := \text{ideal}\{x^{v^+} - x^{v^-} : v \in L\}$ where $v^+, v^- \in \mathbb{N}^d$ denote the nonnegative and nonpositive parts of $v \in \mathbb{N}^d$ with $v = v^+ - v^-$. Assume $L$ is full dimensional with determinant $\det(L) = n$ implying that $I_L$ has length $n$ hence $I_L \in \text{Hilb}_n^d$. An $n$-subset $\lambda = \{\lambda_1, \ldots, \lambda_n\}$ is basic for $I_L$ if and only if it is a set of distinct representatives of the congruence classes of $\mathbb{Z}^d$ modulo $L$, that is, if and only if for all $i \neq j$ we have $\lambda_i - \lambda_j \notin L$: this follows from the fact that $u - v \in L$ if and only if $x^u - x^v \in I_L$ (cf. [11]). Thus, $\Gamma(I_L) = \{\lambda \subset V_n^d : i \neq j \text{ implies } \lambda_i - \lambda_j \notin L\}$; the basis polytope is $\mathcal{B}(I_L) = \text{conv}(\sum \lambda : \lambda \in \Gamma(I_L))$; the state polyhedron is $\mathcal{S}(I_L) = \mathcal{B}(I_L) + \mathbb{R}_+^d$; and the set of initial staircases is $\Lambda(I_L) = \{\lambda \in \Gamma(I_L) : \sum \lambda \text{ vertex of } \mathcal{S}(I_L)\}$. If $\prec$ is a monomial order and $\lambda \in \Lambda(I_L)$ is the initial staircase with $in_\prec(I_L) = I_\lambda$ then for every $u \in \mathbb{N}^d$ we have $[x^u]_\lambda = x^{u_\lambda}$ with $u_\lambda$ the unique representative in $\lambda$ with $u - u_\lambda \in L$. So the reduced Gröbner basis of $I_L$ under $\prec$ is $G_\lambda = \{x^u - x^{u_\lambda} : u \in \min(\bar{\lambda})\}$ and the universal Gröbner basis $\mathcal{U}(I_L) = \cup_{\lambda \in \Lambda(I_L)} G_\lambda$ is efficiently computable by our methods. Here the universal Gröbner basis consists of binomials only, and the set of integer vectors $T := \{u - v : x^u - x^v \in \mathcal{U}(I_L)\}$ which can be read off at once from $\mathcal{U}(I_L)$ is a *universal test set* for the lattice minimization problem that, given any $x \in \mathbb{N}^d$ and any $w \in \mathbb{R}^d_+$, asks for an $x^* \in \mathbb{N}^d$ satisfying $x - x^* \in L$ and minimizing the value $w \cdot x$.

## 6   On the coordinatization of the Hilbert scheme

We conclude with a brief discussion of the embedding of the Hilbert scheme $\text{Hilb}^d_n$ into the Grassmanian of $n$-dimensional subspaces of a vector space of dimension $O(n (\log n)^{d-1})$.

Recall $V^d_n = \cup \binom{\mathbb{N}^d}{n}_{stair}$ is the union of $n$-staircases and $U^d_n = \{v + e_i : v \in V^d_n \text{ and } 0 \leq i \leq d\}$ with $e_0 = 0$ and $e_i$ the $i$-th unit vector in $\mathbb{R}^d$ for $1 \leq i \leq d$. Throughout this section, assume that $U^d_n$ is totally ordered. Let $S := \mathbb{F}[x_1, \ldots, x_d]$, let $S_{U^d_n} := \lin\{x^u : u \in U^d_n\} \subseteq S$ be the $\mathbb{F}$-linear span of the monomials $x^u$ with $u \in U^d_n$, and for each ideal $I \in \text{Hilb}^d_n$ let $I_{U^d_n} := I \cap S_{U^d_n}$.

**Proposition 6.1** *The dimension of the vector subspace $I_{U^d_n}$ of any ideal $I \in \text{Hilb}^d_n$ is $p := |U^d_n| - n$.*

*Proof.* Consider the map $\phi : S_{U^d_n} \to (S/I)$ such that $f = \sum_{u \in U^d_n} a_u x^u \mapsto f + I$. Then $\ker(\phi) = \{f \in S_{U^d_n} : f \in I\} = I_{U^d_n}$. Hence the map $\psi : (S_{U^d_n}/I_{U^d_n}) \to (S/I)$ such that $f + I_{U^d_n} \mapsto f + I$ is injective. Consider any element $g + I$ of $S/I$ and let $[g]_\lambda \in S_{U^d_n}$ be the normal form of $g$ with respect to the reduced Gröbner basis of $I$ corresponding to some initial staircase $\lambda$ of $I$. Then $\psi([g]_\lambda + I_{U^d_n}) = [g]_\lambda + I = g + I$ hence $\psi$ is surjective. Therefore, $S/I$ and $S_{U^d_n}/I_{U^d_n}$ are isomorphic $\mathbb{F}$-spaces hence $\dim(S_{U^d_n}/I_{U^d_n}) = \dim(S/I) = n$ and $\dim(I_{U^d_n}) = |U^d_n| - n = p$. □

**Proposition 6.2** *Any ideal $I \in \text{Hilb}^d_n$ is uniquely determined by its vector subspace $I_{U^d_n}$.*

*Proof.* By definition of $U^d_n$, every reduced Gröbner basis of $I$ lies in $I_{U^d_n}$, implying that the elements of $I_{U^d_n}$ generate $I$ as an ideal. Thus, for any $I, J \in \text{Hilb}^d_n$, $I = J$ if and only if $I_{U^d_n} = J_{U^d_n}$. □

Proposition 6.1 shows that $I_{U^d_n}$ is a $p$-dimensional subspace of $S_{U^d_n}$, i.e., a point on the Grassmanian $\text{Gr}(S_{U^d_n}, p)$ of $p$-dimensional subspaces of the $|U^d_n|$-dimensional vector space $S_{U^d_n}$. This implies that $I_{U^d_n}$, and therefore, by Proposition 6.2, also $I$, inherits standard Plücker coordinates from the Grassmanian, leading to an embedding of the Hilbert scheme in projective space, as follows. Let $(f_i)_{i=1}^p$ be any ordered vector space basis for $I_{U^d_n}$ with $f_i = \sum_{u \in U^d_n} m_{i,u} x^u$, and for a $p$-subset $\nu \subset \mathbb{N}^d$ let $M_\nu := (m_{i,u})_{1 \leq i \leq p, u \in \nu}$ be the corresponding $p \times p$ submatrix of coefficients. The Hilbert scheme is then embedded into projective space by

$$\langle \cdot \rangle : \text{Hilb}^d_n \longrightarrow \mathbb{P}^{\binom{U^d_n}{p}} : I \mapsto \langle I \rangle := \left\{ \det(M_\nu) : \nu \in \binom{U^d_n}{p} \right\} ;$$

it is well known that the Plücker point $\langle I \rangle$ is independent of the choice of basis $(f_i)$ of $I_{U_n^d}$.

We proceed with the dual embedding into the Grassmanian $\mathrm{Gr}(S_{U_n^d}, n)$ of $n$-dimensional subspaces of $S_{U_n^d}$. The space $S_{U_n^d}$ is endowed with the standard monomial basis $\{x^u : u \in U_n^d\}$ and corresponding standard inner product $\langle x^u, x^v \rangle = \delta_{u,v}$. Let $I_{U_n^d}^\perp \subset S_{U_n^d}$ be the orthogonal complement of $I_{U_n^d}$ in $S_{U_n^d}$, which is isomorphic to $S_{U_n^d}/I_{U_n^d}$. Let $(h_i)_{i=1}^n$ be any ordered vector space basis for $I_{U_n^d}^\perp$ with $h_i = \sum_{u \in U_n^d} a_{i,u} x^u$, and for an $n$-subset $\lambda \subset U_n^d$ let $A_\lambda := (a_{i,u})_{1 \leq i \leq n, u \in \lambda}$ be the corresponding $n \times n$ submatrix of coefficients. The dual embedding is then given by

$$\langle \cdot \rangle^\perp \,:\, \mathrm{Hilb}_n^d \longrightarrow \mathbb{P}^{\binom{U_n^d}{n}} \,:\, I \mapsto \langle I \rangle^\perp := \left\{ \det(A_\lambda) \,:\, \lambda \in \binom{U_n^d}{n} \right\}.$$

Next we explain how to actually compute the (dual) Plücker coordinates; note, though, that the number of coordinates is $\Omega\binom{n (\log n)^{d-1}}{n}$ hence exponential, so this computation can *not* be carried out in polynomial time even for $d = 2$. Let $I = \mathrm{ideal}(G_\lambda) \in \mathrm{Hilb}_n^d$ be presented by its Gröbner basis corresponding to some initial staircase $\lambda = \{\lambda_1, \ldots \lambda_n\}$ under some monomial order $\prec$. Compute as in Proposition 4.1 the representative $[x^u]_\lambda = \sum_{i=1}^n a_{i,u} x^{\lambda_i}$ of every $u \in U_n^d$, and for $i = 1, \ldots, n$ let $h_i := \sum_{u \in U_n^d} a_{i,u} x^u$. Then $(h_i)_{i=1}^n$ is an ordered basis of $I_{U_n^d}^\perp$. With $U_n^d$ assumed to be totally ordered, we get the $n \times |U_n^d|$ matrix $A = (a_{i,u})$ of coefficients of this basis; the Plücker coordinates are then read off from the minors of $A$ as $\langle I \rangle^\perp = \{\det(A_\mu) : \mu \in \binom{U_n^d}{n}\}$. In particular, the Plücker coordinate $\langle I \rangle_\mu^\perp := \det(A_\mu)$ is nonzero if and only if $\mu$ is basic for $I$.

For the following classes of ideals, the Plücker coordinates have a natural simple form.

**Monomial Ideals.** Let $I_\lambda = \mathrm{ideal}\{x^v : v \in \min(\bar\lambda)\}$ be the monomial ideal in $\mathrm{Hilb}_n^d$ corresponding to an $n$-staircase $\lambda$. Then (see Section 5) $\lambda$ is the only basic set of $I_\lambda$. Therefore the only nonzero Plücker coordinate is $\langle I \rangle_\lambda^\perp$ and so $\langle I \rangle^\perp = e_\lambda$ is the unit vector in projective space $\mathbb{P}^{\binom{|U_n^d|}{n}}$.

**Point Configurations.** Let $I_C$ be the vanishing ideal of a configuration $C = \{c_1, \ldots, c_n\} \subseteq \mathbb{F}^d$ with $\mathbb{F}$ infinite. For $i = 1, \ldots, n$ let $h_i := \sum_{u \in U_n^d} c_i^u x^u$ with $c_i^u = \prod_{k=1}^d c_{i,k}^{u_k}$. For any polynomial $f = \sum_{u \in U_n^d} m_u x^u \in S_{U_n^d}$, its inner product with $h_i$ satisfies $\langle f, h_i \rangle = \sum_{u \in U_n^d} m_u c_i^u = f(c_i)$, hence is zero if and only if $f$ vanishes on $c_i$. Thus, $f$ is orthogonal to $\mathrm{lin}\{h_1, \ldots, h_n\}$ if and only if it vanishes on $C$, or equivalently, $f \in I_C$. This shows that $(h_i)_{i=1}^n$ is an ordered basis of $(I_C)_{U_d^n}^\perp$.

Let $A = (a_{i,u})$ be the $n \times |U_n^d|$ matrix whose rows are indexed by the points $c_1, \ldots, c_n$ in $C$ and columns by $u \in U_n^d$, with $a_{i,u} := c_i^u$. Then $A$ is the matrix of coefficients of the basis $(h_i)_{i=1}^n$ and hence the Plücker coordinates $\langle I_C \rangle^\perp$ of $I_C$ can be read off from the minors of $A$.

Eric Babson
*University of Washington at Seattle, Seattle, WA 98195, USA.*
*email: babson@math.washington.edu*

Shmuel Onn
*Technion - Israel Institute of Technology, 32000 Haifa, Israel,*
    *and*
*University of California at Davis, Davis, CA 95616, USA.*
*email: onn@ie.technion.ac.il, onn@math.ucdavis.edu,*
    *http://ie.technion.ac.il/∼onn*

Rekha Thomas
*University of Washington at Seattle, Seattle, WA 98195, USA.*
*email: thomas@math.washington.edu*